\documentclass[oneside,12pt,reqno]{amsart}

\usepackage{color}
\usepackage{amsmath,subfigure,latexsym,epsfig,amssymb,amscd,euscript,graphicx,mathrsfs}

\calclayout
\numberwithin{equation}{section}
\newcommand\T[1]{\textrm{#1~}}
\newcommand\cP{\mathcal{P}}
\newcommand\cF{\mathcal{F}}
\newcommand\cE{\mathcal{E}}

\newcommand\cN{\mathcal{N}}
\newcommand\cT{\mathcal{T}}

\newcommand\RE{\mathbb{R}}
\newcommand\KA{\mathbb{K}}

\newcommand\bb[1]{\mathbf{#1}}

\DeclareMathOperator{\Div}{div}
\DeclareMathOperator{\grad}{grad}

\renewcommand{\theenumii}{\arabic{enumii}}
\renewcommand{\labelenumi}{(HG)}
\renewcommand\Gamma{\mathscr T}
\newcommand\nulla[1]{}

\theoremstyle{plain}

\theoremstyle{remark}

\begin{document}

\title[]
{On low order mimetic finite difference methods}
\author{Andrea Cangiani}



\maketitle

\begin{abstract}
These pages review two families of mimetic finite difference methods: the mixed-type methods presented in~\cite{BLSi05} and the nodal methods of~\cite{BBLnodal09}. The purpose of this exercise it to highlight the similitudes underlying the construction of the two families. The comparison prompts the definition of  a piecewise linear postprocessing of the nodal mimetic finite difference solution,  
as it was done for the mixed-type method in~\cite{CM08}.
\end{abstract}

\section{Setting of the problem}
Let $\Omega\subset\RE^3$ be a Lipschitz bounded polyhedral domain.
We consider the following elliptic boundary value problem:
\begin{equation}
\label{eq:strong}
\left\{
 \begin{array}{ll}
 - \nabla\cdot (\KA \nabla u) =g & \mbox{in }\Omega,\\
  u=0 & \mbox{in }\partial \Omega,
 \end{array}
\right.
\end{equation}
where 
$\KA\in\RE^{d\times d}$ is a full symmetric tensor with components 
in $W^{1,{\infty}}(\Omega)$. Strong ellipticity is assumed:
thus there exists two positive constants $\kappa_*$ and $\kappa^*$ such that
\begin{equation}
 \label{eq:ellipticity}
 \kappa_*\|\bb{v}\|^2\le\bb{v}^T\KA(\bb{x})\bb{v}\le\kappa^*\|\bb{v}\|^2\quad\forall\,\bb{v}\in\RE^3,
 \,\mbox{for a.e. }\bb{x}\in\Omega,
\end{equation}
where $\|\cdot\|$ denotes the Euclidean norm of $\RE^3$.

Let $L^2_0(\Omega)$ denote the space consisting of square Lebesgue-integrable functions having zero mean value and
$$H(\Div,\Omega)=\left\{{F}\in (L^2(\Omega))^3:\Div {F}\in L^2(\Omega)\right\}.
$$

In view of the discretization of~\eqref{eq:strong} by both the nodal and mixed-type MFD, we consider the standard and mixed variational formulation of the problem, namely,
\begin{equation}
 \label{eq:weak}
\T{ find $u\in H^1_0(\Omega):$}
   (\KA\nabla u,\nabla v)=(g,v)\quad 
   \forall\,v\in H^1_0(\Omega),
\end{equation}
with $(\cdot,\cdot)$ denting the $L^2$-scalar product, and
\begin{equation}
 \label{eq:mixedweak}
 \begin{tabular}{l}
  \multicolumn{1}{l}{find $(p,F)\in L^2_0(\Omega)\times H(\Div,\Omega)$ s.t.}\\
  $\left\{
  \begin{tabular}{ll}
   $(\KA^{-1}{F},{G}) + (p,\Div {G})=0\quad$ &
   $\forall\,{G}\in H(\Div,\Omega)$,\\
   $(\Div {F},q)=(g,q)$ & $\forall q\in L^2_0(\Omega),$ 
  \end{tabular}
  \right. $
 \end{tabular}
\end{equation}
respectively.

\section{The mimetic finite difference method}
In this section we recall the mimetic finite difference
methods introduced in~\cite{BLSh05,BLSi05} and~\cite{BBLnodal09}.
For more details we refer to the original papers and the references therein.

Let $\cT_h$ be a sequence of non-overlapping conformal decomposition of $\Omega$ into simply-connected polyhedral elements.
For every element $P$ we denote by $|P|$ its
volume and by $h_P$ its diameter. Similarly, for each face $f$ we denote by $|f|$ its area and by $h_f$ its diameter,
and for every edge $e$ we denote by $|e|$ its length.
Depending on the context, $\partial P$ denotes either the boundary of the element $P$ or the union of the element faces. As usual, we set
\[
 h=\max_{P\in \cT_h}h_P.
\]
We assume the following mesh regularity conditions (see \cite{BBLnodal09}).
\begin{enumerate}
 \item {[{\itshape Shape-regularity}]} There exist two positive real numbers $N_s$ and $\rho_s$, independent of $h$, such that every mesh $\cT_h$ 
in the sequence admits a sub-partition $\mathcal{S}_h$ into tetrahedra such that:\medskip
\renewcommand{\labelenumii}{(HG\theenumii)}
\begin{enumerate}
\item every polyhedron $P\in\cT_h$ admits a decomposition $\mathcal{S}_{h|P}$ made of at most $N_s$ tetrahedra;
\item every tetrahedron $T$ of $\mathcal{S}_h$ is {\itshape shape-regular} in the sense that the ratio between the radius 
$r_T$ of the inscribed sphere and its diameter $h_T$ is bounded from below by $\rho_s$, i.e., 
$$
\dfrac{r_T}{h_T}\ge\rho_s>0\quad\forall T\in\mathcal{S}_h.
$$
\end{enumerate}
\renewcommand{\labelenumi}{(ME)}
\item {[{\itshape Star-shaped elements}]} There exists a positive constant $\tau_*$, independent of $h$, 
such that for each element $P$ there exists a point $M_P\in P$ such that $P$ is star-shaped with respect to every point in the ball of center $M_P$ and radius $\tau_*h_P$.
\end{enumerate}

The MFD methods solution is a collection of real values associated to the set of elements, faces, edges, and nodes of the decomposition $\cT_h$.
Following~\cite{BB09}, we thus introduce four discrete spaces $\cP$, $\cF$, $\cE$, and $\cN$
made of collections of real values associated to each element, face, edge, and node, respectively.

The value associated to a face or an edge is to be interpreted as flux or work of vector fields, and thus 
faces and edges are assumed to be given an orientation. For any face $f$ we fix its orientation once and for all by attaching to it a normal vector ${\mathbf{n}}_f$. Further, any edge $e$ with vertices $(V_1,V_2)$, is assumed oriented from $V_1$ to $V_2$.

The notation $\cN_P$ will indicate the restriction of $\cN$ to the nodes belonging to the element $P$, and so on. 
For any element $P$, we also introduce the notation $V_P$ and $f_P$ for the number of vertices and faces, respectively. Further,  we denote by $V_f$ the number of vertices of any face $f$

The nodal MFD method discretizing~\eqref{eq:weak} is built over the set $\cN$ of nodal values 
and the set $\cE$ of edge values. 
The primal discrete differential operator is the gradient
operator $\grad^h:\cN\rightarrow\cE$ defined as follows: for each $u_h\in\cN$, given an edge $e$ with vertices $(V_1,V_2)$,
\[
 (\grad^h u_h)_e=\frac{1}{|e|}(u^e_2-u^e_1),
\]
where $u^e_i$ denotes the value in $u_h$ corresponding to the node $V_i$ of the edge $e$.
This is the discrete gradient operator introduced in~\cite{BBLnodal09},
up to the scaling factor $1/|e|$: here we adopt such scaling for consistency with the definition of the divergence operator below (alternatively, we could have scaled instead the divergence operator, as in~\cite{BB09}).
In order to implement the homogeneous Dirichlet boundary conditions, we shall also need to consider the subset $\cN_0$ of the elements in
$\cN$ that are zero-valued on all vertices belonging to $\partial\Omega$.
Similarly, the elements of the subset $\cE_0$ of $\cE$ are  zero-valued on all the boundary edges.
Notice that the restriction to $\cN_0$ of the operator $\grad^h$ maps into $\cE_0$.

The mixed MFD method for solving~\eqref{eq:mixedweak} is built over the set $\cP$ of elemental values 
and the set $\cF$ of face values. 
The primal discrete differential operator is in this case the divergence operator $\Div^h:\cF\rightarrow \cP$ defined as follows:
for each $F_h\in \cF$, given an element $P$, we set
\[
 (\Div^h F_h)_P=\frac{1}{|P|}\sum_{f\in \partial P} |f| F_f^P,
\]
where $F_f^P=F_f {\mathbf n}_f \cdot{\mathbf n}_f^P$, with ${\mathbf n}_f^P$ the normal of $f$ out of $P$ and $F_f$ representing the value of $F_h$ associated to $f$.

Let us remark that the two operators just defined operate on two different pairs of discrete spaces. In fact, it is possible to define, starting from the gradient and divergence operators above, a {\em derived} divergence and gradient operator, respectively. As these are not needed in the MFD formulation, we refrain to do so.

We shall also need the relevant interpolation operators.
Given any function $p\in L^1(\Omega)$, we define its interpolant $\Pi_\cP p \in \cP$ as
\begin{equation}
  \label{eq:Pint}
 (\Pi_\cP{p})_P=\dfrac{1}{|P|}\int_P p\,\mbox{d}V\quad\T{for all element $P$}.
\end{equation}
For every vector-valued function ${F}\in (L^s(\Omega))^3,\ s>2$, with $\Div {F}\in L^2(\Omega)$, we define its interpolant 
$\Pi_\cF {f}\in \cF$ as
\begin{equation}
  \label{eq:Fint}
 (\Pi_\cF {f})_f=\dfrac{1}{|f|}\int_f {F}\cdot\bb{n}_f\,\mbox{d}S\quad\T{for all face $f$}.
\end{equation}
As the interpolator for $\cE$ shall not be needed, we conclude by
defining the nodal interpolator.
Given any function $u\in H_0^1(\Omega)\cap C^0(\bar{\Omega})$, we define its interpolant $\Pi_\cN u\in \cN$ as
\begin{equation}
 \label{eq:Nint}
	(\Pi_\cN u)_V=u(V)\quad \T{for all vertex $V$}.
\end{equation}

Similarly to finite element methods, an MFD method is defined by restricting a given variational formulation to the MFD discrete set, with the crucial difficulty that the $L^2$-product has to be substituted by {\em discrete} scalar products. Notice that, even when computing a finite element,
the {\em exact} $L^2$-product is substituted by a quadrature formula that is consistent with the degree of accuracy of the method.

The principle used to define MFD methods is indeed that the discrete scalar product has to be consistent, i.e. exact on (the interpolants of) the correct space of polynomials: in the case of low order MFD methods, we shall require exactness on constants. Not surprisingly, the linear and Raviar-Thomas finite element methods are instances of the low order nodal and mixed-type MFD methods, respectively, obtained by a particular choice of MFD scalar products.

In order to derive the methods from the respective consistency conditions, we shall need to work with suitable approximations of the data.
We denote by $\widetilde{g}$  the piecewise constant function
obtained from the forcing function $g$ by averaging over each element $P$ in $\cT_h$. Similarly, $\widetilde{\KA}$ will denote the tensor obtained from $\KA$
by averaging each component over each $P$ in $\cT_h$.

Up to the definition of the relevant scalar products (and linear functionals), we write down the discrete MFD counterparts of the two problems~\eqref{eq:weak}
and~\eqref{eq:mixedweak} as follows.

The nodal MFD method reads:
\begin{equation}
 \label{eq:mfdnodal}
 \T{find} u_h\in\cN_0 \, : \,
   [\grad^h u_h,\grad^h v_h]_{\cE}
=(\tilde{g},v_h)_\cN
	\quad \forall\,{v_h}\in \cN_0.
\end{equation}
Further, we write the mixed MFD method as:
\begin{equation}
\label{eq:mfdmixed}
\begin{tabular}{l}
\multicolumn{1}{l}{find $({F}_h,{p}_h^s)\in \cF\times \cP$:}\medskip\\
$
\left\{
\begin{tabular}{ll}
$[{F}_h,{G}_h]_{\cF} + [{p}_h,\Div^h {G}]_{\cP}=0\quad$ &
$\forall\,{G}_h\in \cF$,\\
${[\Div^h {F}_h,{q}_h]}_{\cP}=-[\Pi_\cP{g},{q}_h]_{\cP}$ & $\forall\,{q}_h\in \cP$. 
\end{tabular}
\right.
$
\end{tabular}
\end{equation}

\section{Scalar products}
The construction of the scalar products is achieved element by element and then summing up the
elemental contributions.
As mentioned earlier, the principle is that we want our scalar products to respect element by element the constants (more precisely, the interpolant of constant functions).

The definition of the product in $\cP$ is straightforward:
\begin{equation}
\label{eq:productq}
[{p},{q}]_{\cP}:=\sum_{E\in\Gamma_h}p_P q_P|P|\quad\forall\,{p},{q}\in \cP.
\end{equation}

Regarding the discrete space $\cN$,  we actually just need to define the linear functional $(\tilde{g},\cdot)_{\cN}$.
To this end, for each $P\in\cT_h$, we introduce the numerical integration formula
\begin{equation}
  \label{eq:quadP}
 \int_P v \, dP\simeq \sum_{i=1}^{V_P} v(V_P^i)\omega_P^i,
\end{equation}
where $\{\omega_P^i\}_{i=1}^{V_P}$ is a set of non-negative weights such that
the quadrature is exact whenever $v$ is a constant.
We then define
\begin{equation}
 \label{eq:linfunN}
(\widetilde{g},v_h)_\cN:=\sum_{P\in\cT_h}\widetilde{g}|_P
 \sum_{i=1}^{V_P} v_h(V_P^i)\omega_P^i.
\end{equation}

We now come to the less obvious problem of the definition of the scalar products mimicking the $H^1_0$-products, namely $[\cdot,\cdot]_{\cF}$
and $[\grad^h \cdot,\grad^h \cdot]_{\cE}$.

The starting point is, in both cases, the Green identity 
\begin{equation}
\label{eq:greenid}
 \int_P 
\Phi\cdot\nabla \phi \, dP=-\int_P\phi \nabla\cdot 
\Phi
\, dP
+\int_{\partial P}
\Phi\cdot  {\mathbf{n}}_P \phi\,dS
\end{equation}
valid for any sufficiently smooth vector-valued function $\Phi$ and scalar function $\phi$.

To deduce consistency conditions for the discrete scalar products,
we specialize the above identity by testing it on the space  of linear polynomials $\mathbb{P}^1(P)$.
To this end, let us consider, as basis of $\mathbb{P}^1(P)$, the set
$\{ b_j\}_{j=0}^3$ given by
\[
\begin{array}{l}
b_0(x)=1,\\
b_{j}(x)=\hat{x}_j\cdot(x-x_P),\quad i=1,\dots,3,
\end{array}
\]
where $\hat{x}_j$
is the $j$-th coordinate vector, and $x_P$ is the barycentre of $P$.

Substituting $\phi=b_0$ in~\eqref{eq:greenid} gives
back the divergence theorem, while with $\phi=b_j$, $j=1,\dots,3$, we obtain
\begin{equation}
\label{eq:conditions_m}
 \int_P 
\widetilde{\KA}^{-1}
\Phi\cdot \widetilde{\KA}\hat{x}_j \, dP=-\int_P  b_j\nabla\cdot 
\Phi
\, dP
+\int_{\partial P}
\Phi\cdot  {\mathbf{n}}_P b_j\,dS,
\end{equation}
where we have also multiplied and divided by $\widetilde{\KA}$ in view of reproducing the weighted $L^2$-product appearing in the first equation of~\eqref{eq:mixedweak}.

This identity suggests to define, for every $G_h\in\cF$, and $j=1,\dots,3$,
\begin{equation}
\label{eq:consistent_mixed}
\begin{array}{ll}
[G_h,\Pi_\cF(\widetilde{\KA}\hat{x}_j)]^P_\cF & := 
\displaystyle{\sum_{f\in\partial P} \int_f G_f^P b_j\, dS}\\
& = \displaystyle{\sum_{f\in\partial P} |f| (x_f-x_P)_j  G_f^P ,}
\end{array}
\end{equation}
with $x_f$ denoting the barycentre of the face $f$. Notice that the volume integral in the right-hand side of~\eqref{eq:conditions_m} disappears due to the fact that $\Div^h G_h$ is constant over $P$.

It easily follows that any scalar product satisfying the so called {\em local consistency} condition~\eqref{eq:consistent_mixed} is exact on the interpolant of constant vectors. Indeed we get that
\begin{equation}\label{eq:exconst_mixed}
[\Pi_\cF(\widetilde{\KA} \hat{x}_i ),\Pi_\cF(\widetilde{\KA} \hat{x}_j) ]^P_\cF
=|P|\widetilde{\KA}_{i,j} %
\\
=\int_P\widetilde{\KA}\nabla b_i\cdot\nabla b_j,
\end{equation}
for all $ i,j=1,\dots,3$.

Substituting $\Phi=\widetilde{\KA}\nabla b_0$ in~\eqref{eq:greenid} gives the triviality $0=0$, while with $\Phi=\widetilde{\KA}\nabla b_j$, $j=1,\dots,3$, we obtain
\begin{equation}
\label{eq:conditions_n}
 \int_P \widetilde{\KA}\nabla b_j\cdot\nabla \phi \, dP
=\int_{\partial P}\widetilde{\KA}\nabla b_j \cdot  {\mathbf{n}}_P \phi\,ds
=\sum_{f\in\partial P}(\widetilde{\KA}_P\hat{x}_j\cdot\mathbf{n}_f^P)\int_f\phi\, ds.
\end{equation}
This time we mimic such identity at the discrete level by requiring that, for every $v_h\in\cN$,
and $j=1,\dots,3$,
\begin{equation}
\label{eq:consistent_nodal}
[\grad^h(v_h),\grad^h(\Pi_\cN b_j)]^P_\cE\equiv [v_h ,\Pi_\cN b_j]^P_\cN
:=\sum_{f\in\partial P}
(\widetilde{\KA}_P\mathbf{n}_f^P)_j
\sum_{l=1}^{V_f} v^f_l\omega^f_l.
\end{equation}
Here, \{$\omega^f_l\}_{l=1}^{V_f}$ represents a set of non-negative weights of a quadrature formula used to approximate the integral over the face $f$. Assuming that such quadrature formula is exact for polynomials of degree $\le 1$, we easily get that
\begin{equation}\label{eq:exconst_nodal}
[\Pi_\cN b_i,\Pi_\cN b_j]^P_\cN
=|P|\widetilde{\KA}_{i,j}=\int_P\widetilde{\KA}\nabla b_i\cdot\nabla b_j
\quad \forall i,j=1,\dots,3.
\end{equation}

{\em Acceptable} MFD scalar products over $\cF$ or $\cN$ are bilinear forms satisfying respectively \eqref{eq:consistent_mixed} or~\eqref{eq:consistent_nodal} which are symmetric and
obey the following scaling properties: there exist two constants $c_*$ and $c^*$ independent
of $P\in\cT_h$ such that
\[
 \begin{array}{l}
\displaystyle{
  c_*\sum_{f\in\partial P} |P|(G^P_f)^2\le [G_h,G_h]^P_\cF\le c^*\sum_{f\in\partial P} |P|(G^P_f)^2
  \quad\forall G_h\in\cF,	     } \\
\displaystyle{
 c_*\sum_{e\in\partial P} |P|(\grad^h v_h)_e^2\le [v_h,v_h]^P_\cN\le c^*\sum_{e\in\partial P} |P|(\grad^h v_h)_e^2
  \quad\forall v_h\in\cN.           }
 \end{array}
\]
The consistency, symmetry and positivity conditions leave some freedom in the definition of the forms, and indeed we shall get a family of MFD scalar products. This is better analyzed by considering the matrices associated to the bilinear forms.

Let $M^P_\cF$ be the $f_P\times f_P$ symmetric matrix related to the form $[\cdot,\cdot]^P_\cF$ and $M^P_\cN$ the $V_P\times V_P$ symmetric matrix related to $[\cdot,\cdot]^P_\cN$.

We translate the conditions~\eqref{eq:consistent_mixed} and~\eqref{eq:consistent_nodal} into algebraic conditions for $M^P_\cF$ and $M^P_\cN$ by introducing the following matrices.
Let $N$ and $R$ be the $f_P\times 3$ matrices   given by
\[
N=[\Pi_\cF(\widetilde{\KA}\hat{x}_1) \, \dots\, \Pi_\cF(\widetilde{\KA}\hat{x}_3)]=
[\mathbf{n}_{f_1}^P \cdots  \mathbf{n}_{f_{f_P}}^P]^T
\widetilde{\KA}_P
\]
(the above equivalence is obtained expressing the definition of the interpolant $\Pi_\cF$) 
and
\[
R=
\left[ 
\begin{array}{c}
|{f_1}|(x_{f_1}-x_P)\\
\vdots \\ 
|{f_P}|(x_{f_{f_P}}-x_P)
\end{array}
\right],
\]
respectively. We write the consistency condition~\eqref{eq:consistent_mixed} as
\begin{equation}
 \label{eq:MNR}
   M^P_\cF N=R.
\end{equation}
Further, let $W$ be the $f_P\times V_P$ matrix collecting on each row the facial quadrature weights
appearing in~\eqref{eq:consistent_nodal} filled with zeros to account for the nodes not belonging to the corresponding face. Then, by introducing the $V_P\times 3$ matrices $A$ and $B$ given by
\[
 A^T=N^TW\quad {\rm and}\quad B=[\Pi_{\cN}b_1 \dots \Pi_{\cN}b_3],
\]
we write the consistency condition~\eqref{eq:consistent_nodal} as
\begin{equation}
 \label{eq:MBA}
   M^P_\cN B=A.
\end{equation}

Let now $C$ be a $f_V\times (f_V-3)$ matrix with columns that span the null space of $N^T$
and $D$ be a $V_P\times (V_P-4)$ matrix with columns that span the null space of
$[\Pi_{\cN}b_0 \cdots \Pi_{\cN}b_3]$. Then the general form of acceptable matrices $M^P_\cF$ and
$M^P_\cN$ is respectively given by
\begin{equation}\label{eq:Mdef}
 M^P_\cF=\frac{1}{|P|}R \widetilde{\KA}_P^{-1}R^T+CU_\cF C^T
\quad {\rm and}\quad 
 M^P_\cN=\frac{1}{|P|}A \widetilde{\KA}_P^{-1}A^T+DU_\cN D^T
\end{equation}
with $U_\cF$ and $U_\cN$ arbitrary symmetric and positive definite matrices of the appropriate scaling. The  dimension of 
$U_\cF$ and $U_\cN$ is $(f_V-3)\times(f_V-3)$ and $(V_P-4)\times (V_P-4)$, respectively.
Notice that, if $P$ is a tetrahedron, than we just get one possible nodal scalar product, namely 
$M^P_\cN=\frac{1}{|P|}A \widetilde{\KA}_P^{-1}A^T$.
As noted in~\cite{BBLnodal09}, when $\cT_h$ is made of tetrahedrons the nodal MFD values coincide with those of the standard
$\mathbb{P}^1$ finite element method (with $\widetilde{\KA}$ used in place of $\KA$ before evaluating the elemental integrals).

\section{Gradient reconstructions}

It is clear that only the first part of the definitions of $M^P_\cF$ and $M^P_\cN$ given in~\eqref{eq:Mdef} acts on the relevant subspace 
of interpolated linear polynomials, and we know that such action is exact.
As our forms operate at the gradient level,
we can use them to define piecewise constant gradient reconstructions which have to be exact
on linear polynomials.
Following~\cite{CM08}, for any vector field $G$ on $P$ we define its reconstruction $G^R$
as
\[
  G^R_i:=\frac{1}{|P|}[\Pi_\cF G,\Pi_\cF(\widetilde{\KA}\hat{x}_i)]^P_\cF\quad\forall i=1,\dots,3,
\]
which is exact on constant fields by~\eqref{eq:exconst_mixed}. Notice that the reconstruction is easily calculated as follows
\begin{equation}
\label{eq:recon_mixed}
 G^R=\frac{1}{|P|}N^T M^P_\cF \left(\Pi_\cF G\right)_P=\frac{1}{|P|} 
\widetilde{\KA}^{-1} R^T  \left(\Pi_\cF G\right)_P.
\end{equation}
The reconstruction formula~\eqref{eq:recon_mixed} was used in~\cite{CM08} to postprocess the mixed MFD method solution $(p_h,F_h)$. 
Indeed, after $(p_h,F_h)$ have been calculated, we can assemble elementwise a piecewise linear second-order accurate solution is given, on each $P\in\cT_h$, by
\[
 p_h^R|_P=p_h|_P+\frac{1}{|P|} \left(\widetilde{\KA}_P^{-1} R^T F_h|_P\right)\cdot(x-x_P). 
\]

Similarly, we can define a reconstructed gradient of any scalar function $v$ over $P$ as
\[
 (\widetilde{\KA}_P\grad^R v)_j:=\frac{1}{|P|}[\Pi_\cN v,\Pi_\cN b_j]^P_\cN, \quad\forall j=1,\dots,3,
\]
which is exact on linear polynomials by~\eqref{eq:exconst_nodal}.
Also in this case we have the following formula for the computation of the reconstructed gradient:
\begin{equation}
\label{eq:recon_nodal}
 \grad^R v=\frac{1}{|P|}\widetilde{\KA}_P^{-1}
B^T M^P_\cN \left(\Pi_\cN v\right)_P 
=\frac{1}{|P|}\widetilde{\KA}_P^{-1}A^T \left(\Pi_\cN v\right)_P.
\end{equation}
This formula may be used to define, starting from the nodal MFD solution $u_h$, a piecewise linear solution $u_h^R$ given, on each $P\in\cT_h$, by
\[
 u_h^R|_P=\frac{1}{|P|}\sum_{i=1}^{V_P}u_h(V_P^i)\omega_P^i+
\frac{1}{|P|}\left(\widetilde{\KA}_P^{-1}A^T u_h|_P\right)\cdot(x-x_P).
\]
If the decomposition $\cT_h$ is made of tetrahedrons,
the exactness of the elemental and facial quadrature formulas on the respective spaces of constant and linear polynomials implies that 
$u_h^R|_P$ is the unique linear polynomial which takes the values $u_h(V_P^i)$ at the elemental vertices. Thus $u_h^R$ coincides 
with the standard $\mathbb{P}^1$ finite element solution. 

It is interesting to rediscover, in the case of tetrahedral elements, the exactness 
of~\eqref{eq:recon_nodal} on linear polynomials. To this purpose, let us consider $\{\varphi_j\}_{j=1}^4$ as the standard basis of 
$\mathbb{P}^1(P)$, so that $\varphi_j(V_P^i)=\delta_{i,j}$. Then, for any $j=1,\dots,4$, denoting by $f_j$ the face opposed to the $j$-th vertex of the tetrahedron, we have
\[
\grad \varphi_j=-\frac{1}{3|P|}|f_j|\mathbf{n}_{f_j}^P.
\]
On the other hand, the reconstructed gradient of $\varphi_j$ is given by
\[
\grad^R \varphi_j=\frac{1}{|P|}\widetilde{\KA}_P^{-1}A^T \left(\Pi_\cN \varphi_j\right)_P
=\frac{1}{|P|}(\widetilde{\KA}_P^{-1}A^T)^j
=\frac{1}{|P|}( [\mathbf{n}_{f_1}^P \cdots  \mathbf{n}_{f_{f_P}}^P]W)^j,
\]
with the symbol $(\cdot)^j$ indicating the $j$-th column of its argument.
Now, the $j$-th column of $W$ is given by $(W)^j=\frac{1}{3}(|f_1| \cdots 0 \cdots|f_{f_P}|)^T$ with the $0$ appearing in the $j$-th position, as these weights  are the only ones ensuring exactness on linear polynomials. It follows that
\[
\grad^R \varphi_j=\frac{1}{|P|}\sum_{\substack{i=1 \\ i\neq j}}^{f_P}
\frac{1}{3}|f_i|\mathbf{n}_{f_i}^P=-\frac{1}{3|P|}|f_j|\mathbf{n}_{f_j}^P.
\]

\section{Discretisation of advective terms: nodal MFD}

We now consider the advection-diffusion problem
\begin{equation}
\label{eq:adve_diff:strong}
\left\{
 \begin{array}{ll}
 - \nabla\cdot (\KA \nabla u) +\beta\cdot \nabla u=g & \mbox{in }\Omega,\\
  u=0 & \mbox{in }\partial \Omega,
 \end{array}
\right.
\end{equation}
where $\beta$ is a given vector field with components in $W^{0,\infty}(\Omega)$.

We want to construct a nodal MFD discretisation of~\eqref{eq:adve_diff:strong}. 
The discretisation of the new, advective, term will be based on
the gradient reconstruction formula~\eqref{eq:recon_nodal} and
on a piecewise constant approximation of the data.

Let $\tilde{\beta}$ be the vector field obtained from $\beta$
by averaging each component over each mesh element $P$. Then
$\tilde{\beta}\cdot \grad^h v_h\in \cP$, for all $v_h\in\cN$.
We define the following nodal MFD method:
\begin{equation}
 \label{eq:adve_diff:mfdnodal}
 \T{find} u_h\in\cN_0 \, : \,
   [\grad^h u_h,\grad^h v_h]_{\cE}
+(\tilde{\beta}\cdot \grad^R u_h,v_h)_{\cN}
=(\tilde{g},v_h)_\cN
	\quad \forall\,{v_h}\in \cN_0.
\end{equation}
We may as well introduce in the method the following streamline-diffusion type stabilizing term:
\begin{equation}
 \label{eq:sd_term}
[\tilde{\tau}_h \tilde{\beta}\cdot\grad^R u_h,\tilde{\beta}\cdot\grad^R v_h]_\cP,
\end{equation}
where $\tilde{\tau}_h$ is a stabilization parameter defined element-wise
in function of the local mesh P\`eclet number.

Notice that, on each element $P$, the reconstructed gradient appearing above gets the following expression:
\begin{equation}\label{eq:rec_grad:formula}
\begin{array}{ll}	
	\vspace{2 mm}
\grad^R v_h|_P&=\frac{1}{|P|}\widetilde{\KA}_P^{-1}A^Tv_h|_P\\
\vspace{2 mm}
&=\frac{1}{|P|}[\mathbf{n}_{f_1}^P \cdots  \mathbf{n}_{f_{f_P}}^P]W 
v_h|_P\\
&=\displaystyle{
\frac{1}{|P|}\sum_{f\in\partial P}{\mathbf{n}}_{f}^P 
\sum_{l=1}^{V_f}v_h(V_f^l)\omega_f^l.
				}
\end{array}
\end{equation}
Thanks to the exactness on linears of the gradient reconstruction, the term~\eqref{eq:sd_term} satisfies a local consistency condition similar 
to~\eqref{eq:consistent_nodal}.
Once again, the starting point is a Green identity. Given a generic smooth function $\phi$ and a linear function 
${\rm p}\in\mathbb{P}^1(P)$, we have
\begin{equation}\label{eq:green:beta}
\int_P (\tilde{\beta}\cdot\nabla {\rm p})(\tilde{\beta}\cdot\nabla\phi)dP
=\int_{\partial P} (\tilde{\beta}\cdot\nabla {\rm p})(\tilde{\beta}\cdot {\mathbf{n}}_P) \phi\,dS
=(\tilde{\beta}\cdot\nabla {\rm p}) \sum_{f\in\partial P} 
(\tilde{\beta}\cdot {\mathbf{n}}_f)\int_{f} \phi\,dS.
\end{equation}
The term~\eqref{eq:sd_term} mimics~\eqref{eq:green:beta} in that, for each
$v_h\in\cN$ and ${\rm p}\in\mathbb{P}^1(P)$,
\begin{equation}\label{eq:beta:local_consistency}
	\begin{array}{ll}
\vspace{2 mm}
[\tilde{\beta}\cdot\grad^R \Pi_\cN {\rm p},\tilde{\beta}\cdot\grad^R v_h ]_\cP^P
&=|P| (\tilde{\beta}\cdot\grad^R \Pi_\cN {\rm p}) (\tilde{\beta}\cdot\grad^R v_h)\\
&= 
\displaystyle{
(\tilde{\beta}\cdot\nabla {\rm p})\sum_{f\in\partial P}(\tilde{\beta}\cdot {\mathbf{n}}_f^P)\sum_{l=1}^{V_f} v(V_f^l) \omega_f^l,
			 }
\end{array}
\end{equation}
{\em cf.} Equation (5.14) in~\cite{BBLnodal09}. 
Indeed~\eqref{eq:beta:local_consistency} follows from~\eqref{eq:rec_grad:formula} and the fact that $\grad^R \Pi_\cN {\rm p}=\nabla {\rm p}$,
and thus $\tilde{\beta}\cdot\grad^R \Pi_\cN {\rm p}=\tilde{\beta}\cdot\nabla {\rm p}$. Actually, this last equality can 
be re-obtained by direct calculation: again from~\eqref{eq:rec_grad:formula} we have that
\[
\begin{array}{ll}	
	\vspace{2 mm}
\tilde{\beta}\cdot\grad^R \Pi_\cN {\rm p}|_P&=
\displaystyle{
\frac{1}{|P|}\sum_{f\in\partial P}(\tilde{\beta}\cdot{\mathbf{n}}_{f}^P) 
\sum_{l=1}^{V_f}{\rm p}(V_f^l)\omega_f^l
				}\\
	\vspace{2 mm}
&=\displaystyle{	\frac{1}{|P|}\sum_{f\in\partial P}\int_f {\rm p}
(\tilde{\beta}\cdot {\mathbf{n}}_{f}^P)\, dS
               }\\
	\vspace{2 mm}
&=\displaystyle{
\frac{1}{|P|}\int_{\partial P}{\rm p}(\tilde{\beta}\cdot {\mathbf{n}}_{f}^P)\, dS
				}\\
    \vspace{2 mm}			
&=\displaystyle{
\frac{1}{|P|}\int_P \tilde{\beta}\cdot\nabla {\rm p}\, dP
=\tilde{\beta}\cdot\nabla {\rm p}.
				}
\end{array}	
\]

 \bibliographystyle{plain}
\bibliography{low}

\begin{thebibliography}{1}

\bibitem{BB09}
F.~Brezzi and A.~Buffa.
\newblock Innovative mimetic discretizations for electromagnetic problems.
\newblock {\em Journal of Computational and Applied Mathematic}, to appear.

\bibitem{BLSh05}
F.~Brezzi, K.~Lipnikov, and M.~Shashkov.
\newblock Convergence of the mimetic finite difference method for diffusion
  problems on polyhedral meshes.
\newblock {\em SIAM J. Numer. Anal.}, 43(5):1872--1896 (electronic), 2005.

\bibitem{BLSi05}
F.~Brezzi, K.~Lipnikov, and V.~Simoncini.
\newblock A family of mimetic finite difference methods on polygonal and
  polyhedral meshes.
\newblock {\em Math. Models Methods Appl. Sci.}, 15(10):1533--1551, 2005.

\bibitem{BBLnodal09}
Franco Brezzi, Annalisa Buffa, and Konstantin Lipnikov.
\newblock Mimetic finite differences for elliptic problems.
\newblock {\em M2AN Math. Model. Numer. Anal.}, 43(2):277--295, 2009.

\bibitem{CM08}
A.~Cangiani and G.~Manzini.
\newblock Flux reconstruction and pressure post-processing in mimetic finite
  difference methods.
\newblock {\em Comp. Meth. Appl. Mech. Engrg.}, 197/9-12:933--945, 2008.

\end{thebibliography}

\end{document}